\documentclass[12pt,a4paper]{article}
\usepackage[utf8]{inputenc}
\usepackage[T1]{fontenc}
\usepackage[english]{babel}
\usepackage{graphicx}
\usepackage{float}
\usepackage{indentfirst}
\usepackage{listings}
\usepackage{caption}
\usepackage{subcaption}
\usepackage{setspace,amsmath,amssymb,amsthm}
\usepackage[left=20mm, top=20mm, right=20mm, bottom=20mm, nohead, footskip=10mm]{geometry}
\usepackage{xcolor}
\usepackage{hyperref}
\definecolor{linkcolor}{HTML}{000080} 
\definecolor{urlcolor}{HTML}{000080}  
\hypersetup{
    pdfstartview=FitH,
    colorlinks=true,      
    citecolor=linkcolor,  
    linkcolor=linkcolor,  
    urlcolor=urlcolor     
}

\newtheorem{theorem}{Theorem}[section]
\newtheorem{lemma}{Lemma}[section]
\newtheorem{definition}{Definition}

\begin{document}

\title{Global Hyperbolicity Criterion for Regular Sub-Lorentzian Manifolds}

\author{Ramazan G.~Beden\thanks{Department of General Control Problems, Faculty of Mechanics and Mathematics, Lomonosov Moscow State University, Moscow, 119991, Russia; email: \texttt{ramazan.gediev@math.msu.ru}}}

\date{} 

\maketitle

\begin{abstract}
In this paper, the global hyperbolicity property of regular sub-Lorentzian manifolds is investigated. A criterion for global hyperbolicity on complete Riemannian manifolds is established in terms of a smooth time function and an admissible growth function. The proof relies on the apparatus of the theory of differential inclusions. The proposed condition relaxes the classical requirements on the differential of a time function. The effectiveness of the obtained criterion is illustrated on classical nonholonomic structures: the Grushin space, the Heisenberg group, and step-2 Carnot groups, for which explicit asymptotic estimates on the expansion rate of causal cones guaranteeing the preservation of global hyperbolicity are found.


\end{abstract}

\section{Introduction}

The concept of global hyperbolicity is one of the most important in the causal theory of spacetimes. Originally introduced by J.~Leray to study the well-posedness of the Cauchy problem for wave equations, it was deeply developed within the framework of classical Lorentzian geometry and general relativity. A fundamental result here belongs to R.~Geroch \cite{geroch1970}, who proved in 1970 that a globally hyperbolic Lorentzian manifold $M$ splits topologically into a direct product $\mathbb{R} \times \Sigma$ (later, A.~Bernal and M.~S\'{a}nchez \cite{bernal2003} strengthened this result to a smooth diffeomorphic splitting).

With the development of geometric control theory and nonholonomic geometry, which received a powerful impetus from the works of M.~Gromov \cite{gromov1996}, a natural need arose to transfer the methods of causal theory to spaces where motions are constrained by a non-integrable subbundle of the tangent bundle. Works on sub-Lorentzian geometry and the study of causality on them belong to M.~Grochowski \cite{grochowski2004}, who investigated the properties of reachable sets and light cones defined on distributions.

In recent years, causality theory has been significantly generalized to the case of arbitrary closed cone structures (see the fundamental survey by E.~Minguzzi \cite{minguzzi2019}). Within this approach, L.\,V.~Lokutsievskiy and A.\,V.~Podobryaev \cite{lokutsievskiy2024} obtained a new existence theorem for longest curves in sub-Lorentzian problems.

Despite a well-developed theoretical foundation, in practice, when testing specific sub-Lorentzian manifolds for global hyperbolicity, a serious problem arises. The criterion (following from Minguzzi's results \cite[Theorem 3.12]{minguzzi2019}) states that a manifold is globally hyperbolic if there exists a smooth time function $T$ whose differential majorizes some norm on the cone, i.e., $d_xT(a) \ge |a|$. However, for nonholonomic structures (such as Carnot groups or Grushin spaces), a direct search for such a global function leads to the necessity of solving systems of partial differential inequalities. As shown in the present paper, these systems often do not have obvious solutions within the class of simple functions.

The goal of this paper is to formulate and prove a new, more flexible criterion for global hyperbolicity for regular sub-Lorentzian manifolds (Theorem \ref{thm:main}), which makes it possible to bypass the aforementioned analytical difficulties.

The idea consists in introducing a growth function $\rho(r)$ depending on the distance to a fixed point. We prove that the classical condition can be relaxed: if there exists a function $T$ satisfying the inequality
\[
d_xT(a) \ge \frac{|a|}{\rho(\operatorname{dist}(x, A))},
\]
and the integral $\int_0^\infty \frac{dr}{\rho(r)}$ diverges, then the sub-Lorentzian manifold retains the property of global hyperbolicity.

\section{Basic Definitions}

In this section, we provide a set of necessary definitions. Let $M$ be a smooth manifold.

\begin{definition}
    Denote by $C_x^+ \subset T_xM$ a nonempty, convex, closed, and pointed cone $(0\in C_x^+$ is the apex of the cone, and the cone contains no lines$)$.
\end{definition}

\begin{definition}
    An antinorm on the cone $C_x^+$ is an upper semicontinuous function $\nu_x:T_xM\rightarrow{\mathbb{R}}\:\sqcup\:\{-\infty\}$ such that:
    \begin{enumerate}
        \item $\nu_{x\vert_{C_x^+}}\geq 0$ and $\nu_x(T_xM\backslash C_x^+)=-\infty$,
        \item $\forall a \in T_xM$ and $\forall \lambda>0$, we have $\nu_x(\lambda a)=\lambda\nu_x(a)$,
        \item $\forall a,b \in T_xM$, we have $\nu_x(a+b)\geq \nu_x(a)+\nu_x(b)$.
    \end{enumerate}
Note that the cone \(C_{x}^{+}\) must be closed due to the upper semicontinuity of the antinorm. Moreover, if \(\nu_{x}|_{C_{x}^{+}} \neq 0\), then \(\nu_{x}|_{\mathrm{ri}C_{x}^{+}} > 0\), where \(\mathrm{ri}C_{x}^{+}\) denotes the relative interior of the cone \(C_{x}^{+}\). Indeed, if a concave function attains its minimum at a point lying in the relative interior of its domain, then this function is constant on its domain (see \cite{lokutsievskiy2024}).
\end{definition}

\begin{definition}
    Suppose that at each point $x\in M$ a cone $C_x^+$ equipped with an antinorm $\nu_x$ is given. The triple $(M,C^+,\nu)$ is called a sub-Lorentzian manifold.
    If the set
    \[
        C^+=\bigcup_{x\in M}C_x^+ \subset TM
    \]
    is closed, and $\nu$ is upper semicontinuous as a function from $TM$ to $\mathbb{R}\: \sqcup \: \{-\infty\}$, then $(M,C^+,\nu)$ is called a regular sub-Lorentzian manifold.
\end{definition}

\begin{definition}
    A curve $y(t):I\rightarrow{M}$ $($where $I\subset\mathbb{R}$ is an interval or a closed segment$)$ is called an admissible path if the curve $y(t)$ is absolutely continuous and $\dot{y}(t)\in C_{y(t)}^+$ for almost all $t\in I$.
\end{definition}

\begin{definition}
    A continuous function $\rho:[0,\infty)\rightarrow{(0,\infty)}$ is called a growth function if
    \[
        \int_0^\infty{\frac{dx}{\rho(x)}}=\infty.
    \]
\end{definition}

\begin{definition}
    A regular sub-Lorentzian manifold $(M, C^+,\nu)$ is called globally hyperbolic if it contains no closed admissible paths and for any compact sets $K_1,K_2\subset M$ the set $R_{K_1}^+\cap R_{K_2}^-$ is compact, where $R^+_{K_1}=\{y:y\in M \text{ and } \exists x_1\in K_1 \text{ and an admissible path from } x_1 \text{ to } y\}$ and $R^-_{K_2}=\{y:y\in M \text{ and } \exists x_2\in K_2 \text{ and an admissible path from } y \text{ to } x_2\}$.
\end{definition}

\section{Global Hyperbolicity Criterion}

\begin{theorem}\label{thm:main}

Let $M$ be a complete Riemannian manifold, and let $(M,C^+,\nu)$ be a regular sub-Lorentzian manifold. Then:

\begin{enumerate}
    \item If there exist a function $T\in C^1(M\rightarrow{\mathbb{R}})$, a growth function $\rho$, and a point $A \in M$ such that
    \begin{equation}
    \forall x\in M\:\forall a\in C_x^+\setminus\{0\} \::\: \tau_x(a)\stackrel{\mathrm{def}}{=}d_xT(a)\geq \frac{|a|}{\rho(\operatorname{dist}(x,A))}, \label{eq:formula}
    \end{equation}
    then $(M,C^+,\nu)$ is globally hyperbolic.
    \item If $(M,C^+,\nu)$ is globally hyperbolic, then there exists a function $T$ such that
    \[
    \forall x\in M\:\forall a\in C_x^+\setminus\{0\} \::\: \tau_x(a)\stackrel{\mathrm{def}}{=}d_xT(a)\geq |a|, \quad \tau_x(a)>\nu_x(a),
    \]
    that is, \eqref{eq:formula} holds with $\rho=1$.
\end{enumerate}
    
\end{theorem}

Let us first prove the first part of the theorem. Consider the set
\[
    U_x=\{a\in C_x^+:\tau_x(a)=1\}.
\]
\begin{lemma}\label{lem:compact}
    The set $U_x\subset T_xM$ is convex and compact for all $x\in M$.
\end{lemma}
\begin{proof}
    Note that 
    \[
        U_x=\{ \xi \in T_xM : \tau_x(\xi)=1 \}\cap C_x^+.
    \] 
    These sets are convex and closed, which means that the set $U_x$ is convex and closed. Let us show that the set $U_x$ is bounded for any $x$. Suppose the contrary: $\exists \xi_n \in U_x : |\xi_n|\rightarrow{\infty}$. Consider 
    \[
        \tilde{\xi}_n=\frac{\xi_n}{|\xi_n|} \in C_x^+;
    \]
    since the cone $C_x^+$ is closed and $\tilde{\xi}_n$ is a bounded sequence, one can pass to a convergent subsequence $\tilde{\xi}_{n_k}\rightarrow{\tilde{\xi}}\neq0$, whence $\tilde{\xi}\in C_x^+$. Furthermore,
    \[
        \tau_x(\tilde{\xi}_{n_k})=\frac{1}{|\xi_{n_k}|}\rightarrow{0},
    \]
    but $\tau_x(\tilde{\xi}_{n_k})\rightarrow{\tau_x(\tilde{\xi})}$. Thus, $\tau_x(\tilde{\xi})=0$, which contradicts the condition 
    \[
        \tau_{x\vert_{C_x^+\backslash\{0\}}}>0.
    \]
\end{proof}

\begin{lemma}\label{lem:semicontinum}
    The set-valued map $x\mapsto U_x$ is upper semicontinuous in $x$.
\end{lemma}
\begin{proof}
    The graph of the set-valued map $x\mapsto U_x$ (as a subset in $TM$) is the intersection of the closed set $C^+$ and the closed set $\{ \xi : \tau_x(\xi)=1 \} \subset TM$, and is therefore closed. Moreover, in a neighborhood of any point of $M$, the sets $U_x$ are bounded by virtue of the inequality $\xi\in U_x:|\xi|\leq \rho(\operatorname{dist}(x,A))$. According to \cite[Ch.~2, \S~5, Lemma~14]{filippov1988}, such a set-valued map is upper semicontinuous.
\end{proof}

\begin{lemma}\label{lem:continious}
    Consider the differential inclusion 
    \[
        \dot{x}\in U_x, \: x(0)=x_0,
    \]
    and let $x(s)$ be its solution for $s\in [0,\tilde{S})$. Then 
    \[
        \exists C:\forall s\in [0,\tilde{S}) \:r(s)=\operatorname{dist}(x(s),A)\leq C.
    \]
\end{lemma}
\begin{proof}
    All subsequent inequalities in the proof of this lemma hold almost everywhere (since admissible paths are absolutely continuous). First,
    \[
        |\dot{x}(s)|\leq \rho(r(s)).
    \]
    Indeed, $\dot{x}(s)\in U_{x(s)}$; consequently, 
    \[
    \tau_{x(s)}(\dot{x}(s))=1 \geq \frac{|\dot{x}(s)|}{\rho(r(s))}
    \]
    by hypothesis. Next,
    \[
        \forall y,z \in M:\operatorname{dist}(y,A)\leq\operatorname{dist}(y,z)+\operatorname{dist}(z,A)\Rightarrow{}\operatorname{dist}(y,A)-\operatorname{dist}(z,A)\leq\operatorname{dist}(y,z).
    \]
    Swapping $z$ and $y$, we obtain
    \[
        \operatorname{dist}(z,A)-\operatorname{dist}(y,A)\leq\operatorname{dist}(z,y)=\operatorname{dist}(y,z) \Rightarrow{} |\operatorname{dist}(y,A)-\operatorname{dist}(z,A)|\leq \operatorname{dist}(y,z).
    \]
    Let us now estimate $r'(s)$. Take $y=x(s+h)$ and $z=x(s)$. Then
    \[
        |r(s+h)-r(s)|\leq \operatorname{dist}(x(s+h),x(s)) \Rightarrow{} \frac{|r(s+h)-r(s)|}{|h|}\leq\frac{\operatorname{dist}(x(s+h),x(s))}{|h|} \Rightarrow{}
    \]
    \[
        \Rightarrow{} |r'(s)|=\lim_{h\rightarrow{0}}\frac{|r(s+h)-r(s)|}{|h|}\leq\lim_{h\rightarrow{0}}\frac{\operatorname{dist}(x(s+h),x(s))}{|h|}=|\dot{x}(s)|\leq \rho(r(s)).
    \]
    From this we get $|r'(s)|\leq \rho(r(s))$. Consider 
    \[
        F(r)=\int_{r_0}^{r}{\frac{du}{\rho(u)}}, \text{ where } r(0)=r_0.
    \]
    The function $F$ is monotonically increasing and tends to infinity as $r\rightarrow{\infty}$ due to the divergence of the integral. Next,
    \[
        \frac{d}{ds}F(r(s))=\frac{r'(s)}{\rho(r(s))}\leq 1\Rightarrow F(r(s))\leq F(r_0)+s.
    \]
    Suppose that $r(s)\rightarrow{\infty}\text{ as } s\rightarrow{\tilde{S}}$. Then $\forall C> r_0,\:\exists s_C<\tilde{S}$ such that $r(s_C)=C$ and $F(C)\leq F(r_0)+s_C\leq F(r_0)+\tilde{S}$, but the left-hand side grows unboundedly as $C$ increases, which contradicts the boundedness of the right-hand side.
\end{proof}

According to Lemma \ref{lem:continious}, no solution of the differential inclusion
\[
    \dot{x}\in U_x, \: x(0)=x_0
\]
escapes to infinity in finite time.

\begin{lemma}\label{lem:reparam}
Let $(M, C^+, \nu)$ be a regular sub-Lorentzian manifold, and let $T \in C^1(M \to \mathbb{R})$ be a time function such that $d_xT(a) > 0$ for all $a \in C^+_x \setminus \{0\}$.

Then for any absolutely continuous admissible curve $x(s)$ defined on a finite interval $s \in [0, S]$ and satisfying the condition $\dot{x}(s) \in C^+_{x(s)}$ almost everywhere, one can introduce a new parameter $t \in [0, t_1]$ such that the reparameterized original curve satisfies the differential inclusion 
\[
    \dot{y}\in U_y.
\]
\end{lemma}
\begin{proof}
Any absolutely continuous curve on a finite interval has finite length and admits an arc-length parameterization (see, for example, \cite[Proposition 2.5.9]{burago2001}). Let $\sigma \in [0, L]$ be the natural parameter (arc length). Denote the reparameterized curve by $\tilde{x}(\sigma)$. It remains absolutely continuous, and its velocity vector has unit length almost everywhere: $|\dot{\tilde{x}}(\sigma)| = 1$. Consequently, it is Lipschitz continuous.

Since the original curve is admissible, $\dot{\tilde{x}}(\sigma) \in C^+_{\tilde{x}(\sigma)}$ holds almost everywhere. Consider in the tangent bundle $TM$ the subset $\mathcal{K} = \{ (x, v) \in C^+ : x \in \tilde{x}([0, L]), \, |v| = 1 \}$ and show that it is compact. Let $(x_n, v_n)$ be an arbitrary sequence in $\mathcal{K}$. Since the image of the curve $\tilde{x}([0, L])$ is compact in $M$, one can extract a convergent subsequence $x_{n_k} \to x \in \tilde{x}([0, L])$. In a local chart near the point $x$, the Riemannian metric is continuous, so the condition $|v_{n_k}| = 1$ implies the boundedness of the sequence $v_{n_k}$. Passing once more to a subsequence, we obtain $(x_{n_k}, v_{n_k}) \to (x, v)$ in $TM$, with $|v| = 1$ by the continuity of the metric. Since in a regular manifold the set $C^+$ is closed in $TM$, the limit $(x, v)$ also belongs to $C^+$, which means $(x, v) \in \mathcal{K}$, proving the compactness of $\mathcal{K}$. The function $(x, v) \mapsto d_xT(v)$, which is continuous on $TM$, is strictly positive on the compact set $\mathcal{K}$ and therefore attains its minimum $c > 0$ on it. Consequently,
\[
\frac{dt}{d\sigma} = d_{\tilde{x}(\sigma)}T(\dot{\tilde{x}}(\sigma)) \ge c > 0 \quad \text{almost everywhere on } [0, L].
\]

Consider the function $t(\sigma) = T(\tilde{x}(\sigma))$. Since its derivative is bounded away from zero by the constant $c$ almost everywhere, the function $t(\sigma)$ is absolutely continuous and strictly increasing. Therefore, there exists an inverse function $\sigma(t)$, which is also strictly monotonic, and its derivative is bounded almost everywhere: $0 \le \frac{d\sigma}{dt} \le \frac{1}{c}$. A function with a bounded derivative satisfies the Lipschitz condition and, as a consequence, is absolutely continuous.

Introduce the curve $y(t) = \tilde{x}(\sigma(t))$. By the chain rule, almost everywhere we have:
\[
\frac{dy}{dt} = \frac{d\tilde{x}}{d\sigma} \cdot \frac{d\sigma}{dt} = \frac{\dot{\tilde{x}}(\sigma)}{dt/d\sigma} = \frac{\dot{\tilde{x}}(\sigma)}{dT(\dot{\tilde{x}}(\sigma))}.
\]
Apply the differential $dT$ to both sides of the equality (using its linearity):
\[
d_{y(t)}T\left(\frac{dy}{dt}\right) = \frac{dT(\dot{\tilde{x}})}{dT(\dot{\tilde{x}})} = 1.
\]
Moreover, since $\dot{\tilde{x}} \in C^+_x$ and $dT(\dot{\tilde{x}}) > 0$, the vector $\frac{dy}{dt}$ is the result of multiplying a vector from the cone by a positive number, and therefore itself belongs to the cone $C^+_x$. 

Thus, $\frac{dy}{dt} \in C^+_{y(t)}$ and $dT\left(\frac{dy}{dt}\right) = \tau_{y(t)}\left(\frac{dy}{dt}\right) = 1$, which by definition means $\frac{dy}{dt} \in U_{y(t)}$.
\end{proof}

Now let us prove the main theorem.
\begin{proof}
    Let us prove item $1$.
    
    First, we show the absence of closed admissible paths. Suppose that $\gamma:[a,b]\rightarrow{M}$ is a closed admissible path $\gamma(a)=\gamma(b)$ that is not identically a single point. Note that 
    \[
        0=T(\gamma(b))-T(\gamma(a))=\int_a^b{\tau_{\gamma(s)}(\dot{\gamma}(s))}\,ds>0,
    \]
    since $\dot{\gamma}(s)\in C^+_{\gamma(s)}$ and $\tau_{\gamma(s)}(\dot{\gamma}(s))>0$ on a set of positive measure. A contradiction is obtained; therefore, there are no closed admissible paths.
    
    For any two compact sets $K_1$ and $K_2$ on the manifold, consider the sets $R_{K_1}^+$ and $R_{K_2}^-$. If $R_{K_1}^+\cap R_{K_2}^-=\varnothing$, then the intersection is compact and there is nothing to prove. If $R_{K_1}^+\cap R_{K_2}^-\neq\varnothing$, then consider:
    \[
        x(s):x(0)=x_1\in K_1, \: x(\tilde{s})=x_2\in K_2,\:\dot{x}(s)\in U_{x(s)}.
    \]
    \[
        \tilde{s}=\int_0^{\tilde{s}}\,ds=\int_0^{\tilde{s}}\tau_{x(s)}(\dot{x}(s))\,ds=\int_0^{\tilde{s}}\,d(T(x(s)))=T(x_2)-T(x_1)\leq \max_{y\in K_2}T(y)-\min_{y\in K_1}T(y)\stackrel{\mathrm{def}}{=}s_1,
    \]
    whence we obtain that any admissible curve satisfying $\dot{x}(s)\in U_{x(s)}$ goes from $K_1$ to $K_2$ in time $\tilde{s}\leq s_1$. By Lemma \ref{lem:reparam}, any other admissible curve can be reparameterized so as to satisfy the differential inclusion $\dot{x}\in U_x$. By Filippov's theorem for smooth manifolds, the sets
    \[
        R^+_{K_1}(s_1)=\bigcup_{\{x(s):\dot{x}(s)\in U_{x(s)},\:x(0)\in K_1,\:s\in[0,s_1]\}}x(s) \subset M,
    \]
    \[
        R^-_{K_2}(s_1)=\bigcup_{\{x(s):\dot{x}(s)\in -U_{x(s)},\:x(0)\in K_2,\:s\in[0,s_1]\}}x(s) \subset M
    \]
    are compact. Then their intersection $R^+_{K_1}(s_1)\cap R^-_{K_2}(s_1) \subset M$ is also compact. Let us show that
    \[
        R_{K_1}^+\cap R_{K_2}^-=R^+_{K_1}(s_1)\cap R^-_{K_2}(s_1).
    \]
    Let $x\in R_{K_1}^+\cap R_{K_2}^-$; then there exist admissible paths $\gamma_1$ from some $x_1\in K_1$ to $x$ and $\gamma_2$ from $x$ to some $x_2\in K_2$. The curve $\gamma=\gamma_1\cup \gamma_2$ is admissible and goes from $x_1$ to $x_2$; it can be reparameterized according to Lemma \ref{lem:reparam}, and we obtain that the point $x$ lies in $R^+_{K_1}(s_1)\cap R^-_{K_2}(s_1)$.
    Conversely, let $x\in R^+_{K_1}(s_1)\cap R^-_{K_2}(s_1)$; then $\exists x_+(s):x_+(0)\in K_1,\: x_+(\hat{s}_+)=x$ ($\hat{s}_+\in[0,s_1]$) and $\exists x_-(s):x_-(0)\in K_2,\: x_-(\hat{s}_-)=x$ ($\hat{s}_-\in[0,s_1]$), whence it follows that $x\in R_{K_1}^+\cap R_{K_2}^-$. We conclude that the set $R_{K_1}^+\cap R_{K_2}^-$ is compact.

    Item 2 follows from \cite[Theorem 3.12]{minguzzi2019}, where the constant function equal to $1$ is taken as the function $\rho$.
\end{proof}

\section{Examples}
\subsection{Grushin Space}

Consider the three-dimensional spacetime $M = \mathbb{R}_t \times \mathbb{R}^2_{x,y}$. Define the Grushin distribution generated by the vector fields:
\[
X_0 = \frac{\partial}{\partial t}, \quad X_1 = \frac{\partial}{\partial x}, \quad X_2 = x \frac{\partial}{\partial y}.
\]
The distribution $\Delta = \text{span}\{X_0, X_1, X_2\}$ has a singularity on the plane $x=0$, where its rank drops to two. Any velocity vector in coordinates $(t,x,y)$ has the form $a = (v_t, v_x, v_y)$. If $a \in \Delta$, then $a = (v_t, u_1, xu_2)$.

Define on $M$ the Riemannian metric $|a|^2 = v_t^2 + v_x^2 + v_y^2 = v_t^2 + u_1^2 + u_2^2 x^2$. 
Let a smooth strictly positive function $\omega(x,y)$ responsible for the expansion of light cones be given. Define the causal cone with respect to the sub-Riemannian norm of the spatial part:
\[
C^+_{(t,x,y)} = \left\{ a = v_t X_0 + u_1 X_1 + u_2 X_2 \ : \; v_t \ge \frac{1}{\omega(x,y)} \sqrt{u_1^2 + u_2^2} \right\}.
\]

Equip the cone with a strictly positive (inside the cone) antinorm:
\[
\nu(a) = v_t - \frac{1}{\omega(x,y)} \sqrt{u_1^2 + u_2^2} \quad \text{for } a \in C^+,
\]
and $\nu(a) = -\infty$ outside $C^+$. Note that on the singular plane $x=0$, where $X_2 = 0$, the sub-Riemannian norm on $\Delta$ is defined via the infimum of $\sqrt{u_1^2 + u_2^2}$ over all representations $a = v_t X_0 + u_1 X_1 + u_2 X_2$ (which corresponds to setting $u_2 = 0$ at $x=0$), so that $C^+$ coincides with the closure of $C^+|_{x \neq 0}$ in $TM$ and $\nu$ is upper semicontinuous on $TM$.

Direct application of the classical criterion (Condition 2 of Theorem \ref{thm:main}) requires finding a global function $s(t,x,y)$ satisfying the inequality $ds(a) \ge |a|$.

Let us try to solve the following inequality:
\begin{equation}
    ds(a) \ge |a|.
    \label{eq:main_inequality}
\end{equation}

Write the left-hand and right-hand sides of inequality \eqref{eq:main_inequality} in the coordinates of the vector $a = v_t X_0 + u_1 X_1 + u_2 X_2$. The full Riemannian length of the vector is $|a| = \sqrt{v_t^2 + u_1^2 + x^2 u_2^2}$, and the differential of the unknown function $s$ on this vector has the form $ds(a) = s_t v_t + s_x u_1 + s_y (x u_2)$. Then the inequality can be rewritten as:
\begin{equation}
    s_t v_t + s_x u_1 + x s_y u_2 \ge \sqrt{v_t^2 + u_1^2 + x^2 u_2^2}.
    \label{eq:expanded_inequality}
\end{equation}

This condition must hold for all vectors inside and on the boundary of the cone $C^+$. Consider vectors lying strictly on the boundary of the cone, where $v_t$ takes its minimum possible value:
\[
    v_t = \frac{1}{\omega(x,y)} \sqrt{u_1^2 + u_2^2}.
\]
Pass to polar coordinates by setting $u_1 = \omega(x,y) v_t \cos \theta$ and $u_2 = \omega(x,y) v_t \sin \theta$, where $\theta \in [0, 2\pi]$. Substitute these expressions into inequality \eqref{eq:expanded_inequality}:
\begin{equation*}
    v_t \Big( s_t + \omega s_x \cos \theta + x \omega s_y \sin \theta \Big) \ge \sqrt{v_t^2 + \omega^2 v_t^2 \cos^2 \theta + x^2 \omega^2 v_t^2 \sin^2 \theta}.
\end{equation*}

Since $v_t > 0$, we can divide the inequality by $v_t$. As a result, we obtain a partial differential inequality that must hold for any angle $\theta$:
\begin{equation}
    s_t + \omega (s_x \cos \theta + x s_y \sin \theta) \ge \sqrt{1 + \omega^2 (\cos^2 \theta + x^2 \sin^2 \theta)}.
    \label{eq:theta_inequality}
\end{equation}

Consider two limiting cases of the spatial velocity direction.

Case 1: Motion along the non-degenerate field $X_1$ ($\theta = 0$ and $\theta = \pi$).
Substituting these angles, we obtain the system:
\begin{equation*}
    \begin{cases}
        s_t + \omega s_x \ge \sqrt{1 + \omega^2}, \\
        s_t - \omega s_x \ge \sqrt{1 + \omega^2}
    \end{cases}
    \implies \quad s_t - \omega(x,y) |s_x| \ge \sqrt{1 + \omega^2(x,y)}.
\end{equation*}

Case 2: Motion along the degenerate field $X_2$ ($\theta = \pi/2$ and $\theta = 3\pi/2$).
A similar substitution yields:
\begin{equation*}
    \begin{cases}
        s_t + x \omega s_y \ge \sqrt{1 + x^2 \omega^2}, \\
        s_t - x \omega s_y \ge \sqrt{1 + x^2 \omega^2}
    \end{cases}
    \implies \quad s_t - \omega(x,y) |x s_y| \ge \sqrt{1 + x^2 \omega^2(x,y)}.
\end{equation*}
From this follows the condition on the time derivative:
\begin{equation}
    s_t \ge \sqrt{1 + x^2 \omega^2(x,y)} + \omega(x,y) |x s_y|.
    \label{eq:st_bound}
\end{equation}

Let us analyze the obtained constraint \eqref{eq:st_bound} in the case of interest where the cones do not shrink at infinity, i.e., $\omega(x,y) \ge c > 0$ (or, more generally, when $|x|\omega(x,y) \to \infty$ as $|x| \to \infty$):
\begin{enumerate}
    \item Let us seek a solution in the form $s(t,x,y) = C t + f(x,y)$, where $C = \text{const}$ (in particular, $s(t,x,y) = t$ for $C=1$, $f \equiv 0$). Then $s_t = C$. However, for $\omega(x,y) \ge c > 0$ (and more generally whenever $x\omega(x,y) \to \infty$ as $|x| \to \infty$), the right-hand side of inequality \eqref{eq:st_bound} contains the term $\sqrt{1 + x^2 \omega^2(x,y)}$, which grows unboundedly as $|x| \to \infty$. Consequently, the constant $C$ cannot majorize the right-hand side globally on all of $\mathbb{R}^2$.
    \item Let $s(t,x,y) = \alpha(x,y) t + f(x,y)$. Then $s_x = \alpha_x t + f_x$. Substituting this into Case 1, we see that for $\alpha_x \neq 0$ the term $\omega |s_x| = \omega |\alpha_x t + f_x|$ forces the left-hand side $s_t - \omega |s_x|$ to tend to $-\infty$ as $t \to \pm\infty$. Similarly, from Case 2 for $x \neq 0$ and $\alpha_y \neq 0$, the left-hand side $s_t - \omega |x s_y|$ tends to $-\infty$ as $t \to \pm\infty$. For the inequalities to hold for all $t \in \mathbb{R}$, we must set $\alpha_x = 0$ and $\alpha_y = 0$, i.e., $\alpha = \text{const}$, which brings us back to item 1.
\end{enumerate}
Thus, for $\omega(x,y) \ge c > 0$, a direct search in the simplest classes of functions ($t$, $Ct + f(x,y)$, $\alpha(x,y)t + f(x,y)$) does not yield a global solution to inequality \eqref{eq:main_inequality}.

Let us show that Condition 1 of Theorem \ref{thm:main} allows us to avoid solving this inequality (moreover, for an arbitrary smooth function $\omega(x,y) > 0$). Choose the time function $T(t,x,y) = t$; then $dT(a) = v_t$. For any vector $a \in C^+_{(t,x,y)}$, let us estimate its full Riemannian length:
\begin{align*}
|a|^2 &= v_t^2 + u_1^2 + u_2^2 x^2 \\
      &\le v_t^2 + \max(1, x^2) (u_1^2 + u_2^2) \\
      &\le v_t^2 + \max(1, x^2) \omega^2(x,y) v_t^2 \\
      &= v_t^2 \Big( 1 + \max(1, x^2)\omega^2(x,y) \Big).
\end{align*}

Extracting the square root, we obtain an estimate for the differential of time:
\[
dT(a) = v_t \ge \frac{|a|}{\sqrt{1 + \max(1, x^2)\omega^2(x,y)}}.
\]

Denote by $r = \sqrt{x^2+y^2}$ the distance from the point $(x,y)$ to the origin in the spatial plane; then $|x| \le r$, and the Riemannian distance in $M$ satisfies the estimate $\operatorname{dist}((t,x,y), 0) \ge r$. Set $\bar{\omega}(r) = \sup_{x^2+y^2 \le r^2} \omega(x,y)$. 

Then we can define an admissible growth function in the form:
\[
\rho(r) = \sqrt{1 + \max(1, r^2)\bar{\omega}^2(r)}.
\]
Since the function $\rho(r)$ is continuous and non-decreasing, we have $\rho(r) \le \rho(\operatorname{dist}((t,x,y), 0))$, whence $dT(a) \ge |a|/\rho(\operatorname{dist}((t,x,y), 0))$. According to Condition 1 of Theorem \ref{thm:main}, if the integral diverges:
\[
\int_0^\infty \frac{dr}{\sqrt{1 + \max(1, r^2)\bar{\omega}^2(r)}} = \infty,
\]
then this sub-Lorentzian Grushin space is globally hyperbolic. 

Asymptotically as $r \to \infty$, the function behaves as $\rho(r) \sim r \bar{\omega}(r)$. Consequently, global hyperbolicity is preserved, for example, under a logarithmic expansion of the cones $\bar{\omega}(r) \sim \ln(r)$, since $\int (r \ln r)^{-1} dr = \infty$; however, for any power-law expansion $\bar{\omega}(r) \sim r^\varepsilon$ ($\varepsilon > 0$), this choice of time function does not provide an answer regarding the preservation of global hyperbolicity.

\subsection{Heisenberg Group}

Consider the four-dimensional manifold $M = \mathbb{R}_\tau \times \mathbb{H}^1$, where $\mathbb{H}^1$ is the three-dimensional Heisenberg group. Introduce coordinates $(\tau, x, y, z)$ and the sub-Lorentzian distribution $\Delta = \text{span}\{X_0, X_1, X_2\}$ generated by the fields:
\[
X_0 = \frac{\partial}{\partial \tau}, \quad X_1 = \frac{\partial}{\partial x} - \frac{y}{2}\frac{\partial}{\partial z}, \quad X_2 = \frac{\partial}{\partial y} + \frac{x}{2}\frac{\partial}{\partial z}.
\]
Any admissible velocity vector $a \in \Delta$ has the form $a = v_\tau X_0 + u_1 X_1 + u_2 X_2$. In standard coordinates, its components are:
\[
a = \left( v_\tau, \; u_1, \; u_2, \; \frac{1}{2}(u_2 x - u_1 y) \right).
\]
Define the Euclidean metric on $M$ as $|a|^2 = v_\tau^2 + v_x^2 + v_y^2 + v_z^2$. Let a smooth strictly positive function $\omega(x,y,z)$ be responsible for the expansion of the cones. Define the causal cone:
\[
C^+_{(\tau, x,y,z)} = \left\{ a \in \Delta \ : \ v_\tau \ge \frac{1}{\omega(x,y,z)} \sqrt{u_1^2 + u_2^2} \right\}.
\]

Similarly to the previous example, define the antinorm on the cone:
\[
\nu(a) = v_\tau - \frac{1}{\omega(x,y,z)} \sqrt{u_1^2 + u_2^2} \quad \text{for } a \in C^+,
\]
setting it equal to $-\infty$ outside the cone.

Classical verification of global hyperbolicity (Condition 2 of Theorem \ref{thm:main}) encounters significant difficulties due to the nonholonomic structure of the fields $X_1, X_2$ (in particular, $[X_1, X_2] = \partial_z$). 

Let us apply Condition 1 of Theorem \ref{thm:main}. Choose the time function $T(\tau, x, y, z) = \tau$, whence $dT(a) = v_\tau$. Let us estimate the Riemannian length of a vector $a \in C^+$. 
For the component $v_z$, apply the Cauchy--Bunyakovsky--Schwarz inequality:
\[
v_z^2 = \frac{1}{4}(u_2 x - u_1 y)^2 \le \frac{1}{4}(u_1^2 + u_2^2)(x^2 + y^2).
\]
Then the square of the full length is estimated as:
\begin{align*}
|a|^2 &= v_\tau^2 + u_1^2 + u_2^2 + v_z^2 \\
      &\le v_\tau^2 + (u_1^2 + u_2^2) + \frac{1}{4}(u_1^2 + u_2^2)(x^2 + y^2) \\
      &= v_\tau^2 + (u_1^2 + u_2^2)\left( 1 + \frac{x^2+y^2}{4} \right).
\end{align*}
Using the definition of the cone $u_1^2 + u_2^2 \le \omega^2 v_\tau^2$, we obtain:
\[
|a|^2 \le v_\tau^2 \left[ 1 + \omega^2(x,y,z) \left( 1 + \frac{x^2+y^2}{4} \right) \right].
\]
Extracting the square root, we find an estimate for the differential of time:
\[
dT(a) = v_\tau \ge \frac{|a|}{\sqrt{1 + \omega^2(x,y,z) \left( 1 + \frac{x^2+y^2}{4} \right)}}.
\]

Let $r = \sqrt{x^2+y^2+z^2}$ be the distance in space. Note that $x^2+y^2 \le r^2$ and $\operatorname{dist}((\tau,x,y,z), 0) \ge r$. Denote $\bar{\omega}(r) = \sup_{x^2+y^2+z^2\le r^2} \omega(x,y,z)$. 
Then the admissible non-decreasing growth function has the form:
\[
\rho(r) = \sqrt{1 + \bar{\omega}^2(r) \left( 1 + \frac{r^2}{4} \right)}.
\]
According to Theorem \ref{thm:main}, this space is globally hyperbolic if the integral diverges:
\[
\int_0^\infty \frac{dr}{\sqrt{1 + \bar{\omega}^2(r) \left( 1 + \frac{r^2}{4} \right)}} = \infty.
\]
Global hyperbolicity is preserved for constant cones or under their logarithmic expansion; however, if the expansion function $\omega(x,y,z)$ has power-law growth (so that the integral converges), then this choice of time function yields no information about global hyperbolicity.

\subsection{Step-2 Carnot Group}

Let $\mathbb{G}$ be a connected, simply connected nilpotent Lie group whose Lie algebra admits a direct sum decomposition $\mathfrak{g} = V_1 \oplus V_2$, where $[V_1, V_1] = V_2$ and $[V_1, V_2] = [V_2, V_2] = 0$. Let $m = \dim(V_1)$ and $k = \dim(V_2)$.

Consider the ($1+m+k$)-dimensional manifold (spacetime) $M = \mathbb{R}_\tau \times \mathbb{G}$. The sub-Lorentzian distribution $\Delta$ on $M$ is generated by the time vector field $X_0 = \partial_\tau$ and the horizontal fields of the group $X_1, \dots, X_m$, which in standard coordinates take the form:
\[
X_i = \frac{\partial}{\partial x_i} + \frac{1}{2} \sum_{j=1}^m \sum_{l=1}^k B_{ij}^l x_j \frac{\partial}{\partial z_l}, \quad i=1, \dots, m,
\]
where $B^l$ are skew-symmetric matrices of size $m \times m$ (the commutators of the basis fields). 

Any velocity vector $a \in \Delta$ is written as $a = v_\tau X_0 + \sum_{i=1}^m u_i X_i$. Denote the horizontal velocity vector by $u \in \mathbb{R}^m$. Note that the $l$-th velocity component along the fibers of $V_2$ (the axes $z_l$) is computed as:
\[
v_{z_l} = \frac{1}{2} \sum_{i=1}^m u_i \sum_{j=1}^m B_{ij}^l x_j = \frac{1}{2} \langle u, B^l x \rangle_{\mathbb{R}^m}.
\]

Define the Euclidean Riemannian metric on $M$ by the formula $|a|^2 = v_\tau^2 + |v_x|^2 + |v_z|^2$ (for $a \in \Delta$ we have $v_x = u$, so that $|a|^2 = v_\tau^2 + |u|^2 + |v_z|^2$). Let $\omega(x,z) > 0$ be a smooth cone expansion function. Define the cone:
\[
C^+_{(\tau, x,z)} = \left\{ a \in \Delta \ : \  v_\tau \ge \frac{1}{\omega(x,z)} |u| \right\},
\]
equipped with the antinorm $\nu(a) = v_\tau - \frac{|u|}{\omega(x,z)}$ inside $C^+$ and $\nu = -\infty$ outside $C^+$.

Let us apply Condition 1 of Theorem \ref{thm:main}. Set $T(\tau, x, z) = \tau$, whence $dT(a) = v_\tau$. Let us estimate $|a|$. Using the Cauchy--Bunyakovsky--Schwarz inequality and the properties of the operator norm of matrices $\|B^l\|$, we obtain an estimate for $v_{z_l}$:
\[
v_{z_l}^2 = \frac{1}{4} \langle u, B^l x \rangle^2 \le \frac{1}{4} |u|^2 |B^l x|^2 \le \frac{1}{4} \|B^l\|^2 |u|^2 |x|^2.
\]
Summing over all $l = 1, \dots, k$, we find the square of the velocity along the vertical directions:
\[
|v_z|^2 = \sum_{l=1}^k v_{z_l}^2 \le C^2 |u|^2 |x|^2,
\]
where $C^2 = \frac{1}{4}\sum_{l=1}^k \|B^l\|^2$ is a constant depending solely on the Lie algebra of the group $\mathbb{G}$.

Let us estimate the square of the Riemannian length of a vector $a \in C^+$:
\begin{align*}
|a|^2 &= v_\tau^2 + |u|^2 + |v_z|^2 \\
      &\le v_\tau^2 + |u|^2 + C^2 |u|^2 |x|^2 \\
      &= v_\tau^2 + |u|^2 \left( 1 + C^2 |x|^2 \right).
\end{align*}
Taking into account the cone inequality $|u|^2 \le \omega^2 v_\tau^2$, we have:
\[
|a|^2 \le v_\tau^2 \Big( 1 + \omega^2(x,z) (1 + C^2 |x|^2) \Big).
\]
Extracting the square root, we obtain a strict analytical estimate for the differential of time:
\[
dT(a) = v_\tau \ge \frac{|a|}{\sqrt{1 + \omega^2(x,z) (1 + C^2 |x|^2)}}.
\]

Let $r = \sqrt{|x|^2+|z|^2}$ be the Euclidean distance in space, so that $|x| \le r$ and $\operatorname{dist}((\tau,x,z), 0) \ge r$. Taking the supremum of the expansion function $\bar{\omega}(r) = \sup_{|x|^2+|z|^2 \le r^2} \omega(x,z)$, we obtain the desired non-decreasing growth function:
\[
\rho(r) = \sqrt{1 + \bar{\omega}^2(r) (1 + C^2 r^2)}.
\]

According to Theorem \ref{thm:main}, the original space is globally hyperbolic if
\[
\int_0^\infty \frac{dr}{\sqrt{1 + \bar{\omega}^2(r)(1 + C^2 r^2)}} = \infty.
\]

Otherwise, the time function chosen in this way does not provide an answer to the question of global hyperbolicity.

\section{Conclusion}

In the present paper, a new criterion for the global hyperbolicity of regular sub-Lorentzian manifolds is proposed. The introduction of an admissible growth function made it possible to relax the classical requirements on a global time function.

The main advantage of the proved theorem lies in its direct practical applicability to nonholonomic structures. The classical approach, which requires constructing a time function whose differential majorizes the norm of the velocity vector, leads in the sub-Lorentzian case to systems of partial differential inequalities.

Thanks to the new criterion, this analytical difficulty has been bypassed. Using the examples of the Grushin space, the Heisenberg group, and arbitrary step-2 Carnot groups, it is demonstrated that the problem of testing for global hyperbolicity now reduces to algebraic estimates of the asymptotic behavior of causal cones at infinity. Thus, the mathematical apparatus developed in this paper provides a convenient tool for studying the causal structure of a wide class of spaces.

\section*{Acknowledgments}
The author expresses his deep gratitude to his scientific advisor Professor L.\,V.~Lokutsievskiy for posing the problem, fruitful discussions, and constant attention to this work.

\end{document}